\documentclass[11pt]{article}
\usepackage{amsfonts}
\usepackage{mathrsfs}
\usepackage[dvips]{graphics}
\usepackage[dvips]{color}
\usepackage{amsthm}
\usepackage[]{amsmath}
\usepackage{CJK}
\usepackage{indentfirst}
\newtheorem{proposition}{Proposition}[section]

\begin{document}
\begin{CJK*}{GBK}{song}
\CJKindent

\centerline{\textbf{\LARGE{The Square Trees in the Tribonacci Sequence}}}

\vspace{0.2cm}

\centerline{Huang Yuke\footnote[1]{School of Mathematics and Systems Science, Beihang University (BUAA), Beijing, 100191, P. R. China. E-mail address: huangyuke07@tsinghua.org.cn,~hyg03ster@163.com (Corresponding author).}
~~Wen Zhiying\footnote[2]{Department of Mathematical Sciences, Tsinghua University, Beijing, 100084, P. R. China. E-mail address: wenzy@tsinghua.edu.cn.}}

\vspace{1cm}

\noindent\textbf{Abstract:} 
The Tribonacci sequence $\mathbb{T}$ is the fixed point of the substitution $\sigma(a,b,c)=(ab,ac,a)$.
In this note, we get the explicit expressions of all squares, and then establish the tree structure of the positions of repeated squares in $\mathbb{T}$, called square trees.
Using the square trees,
we give a fast algorithm for counting the number of repeated squares in $\mathbb{T}[1,n]$ for all $n$, where $\mathbb{T}[1,n]$ is the prefix of $\mathbb{T}$ of length $n$. Moreover we get explicit expressions for some special $n$ such as $n=t_m$ (the Tribonacci number) etc., which including some known results such as H.Mousavi and J.Shallit\cite{MS2014}.

\vspace{0.2cm}


\vspace{0.2cm}

\noindent\textbf{2010 MR Subject Classification:} 11B85; 68Q45

\section{Introduction}

Let $\mathcal{A}=\{a,b,c\}$ be a three-letter alphabet.
The Tribonacci sequence $\mathbb{T}$ is the fixed point beginning with $a$ of the substitution $\sigma(a,b,c)=(ab,ac,a)$.
As a natural generalization of the Fibonacci sequence, $\mathbb{T}$ has
been studied extensively by many authors, see \cite{DHS2013,MS2014,RSZ2010,TW2007}.

Let $\omega$ be a factor of $\mathbb{T}$, denoted by $\omega\prec\mathbb{T}$.
Let $\omega_p$ be the $p$-th occurrence of $\omega$.
If the factor $\omega$ and integer $p$ such that $\omega_p\omega_{p+1}\prec\mathbb{T}$, we call $\omega_p\omega_{p+1}$ a square  of $\mathbb{T}$.
As we know, $\mathbb{T}$ contains no fourth powers. The properties of squares and cubes are objects of a great interest in mathematics and computer science etc.

We denote by $|\omega|$ the length of $\omega$. Denote $|\omega|_\alpha$ the number of letter $\alpha$ in $\omega$, where $\alpha\in\mathcal{A}$.
Let $\tau=x_1\cdots x_n$ be a finite word (or $\tau=x_1x_2\cdots$ be a sequence).
For any $i\leq j\leq n$, define $\tau[i,j]=x_ix_{i+1}\cdots x_{j-1}x_j$.
By convenient, denote $\tau[i]=\tau[i,i]=x_i$, $\tau[i,i-1]=\varepsilon$ (empty word).
Denote $T_m=\sigma^m(a)$ for $m\geq0$, $T_{-2}=\varepsilon$, $T_{-1}=c$.
Denote $t_m=|T_m|$ for $m\geq-2$, called the $m$-th Tribonacci number.
Denote by $\delta_m$ the last letter of $T_m$ for $m\geq-1$.

Denote $A(n)=\sharp\{\omega,p):\omega_p\omega_{p+1}\prec\mathbb{T}[1,n]\}$ the number of repeated squares in $\mathbb{T}[1,n]$.
In 2014, H.Mousavi and J.Shallit\cite{MS2014} gave expression of $A(t_m)$, which they proved by mechanical way. In \cite{HW2016-2}, we give a fast algorithm for counting the number of repeated squares in each prefix of the Fibonacci sequence.
In this note, we give a fast algorithm for counting $A(n)$ for all $n$.
In Section 2, we establish the tree structure of the positions of repeated squares in $\mathbb{T}$, called the square trees.
Section 3 is devoted to  give a fast algorithm for counting $A(n)$. As a special case, we get expression of $A(t_m)$ in Section 4.

The main tools of the paper are ``kernel word'' and ``gap sequence", which introduced and studied in \cite{HW2015-2}.
We define the kernel numbers that $k_{0}=0$, $k_{1}=k_{2}=1$, $k_m=k_{m-1}+k_{m-2}+k_{m-3}-1$ for $m\geq3$.
The kernel word with order $m$ is defined as
$K_1=a$, $K_2=b$, $K_3=c$, $K_m=\delta_{m-1}T_{m-3}[1,k_m-1]$ for $m\geq4$.
Using the property of gap sequence, we can determine the positions of all $\omega_p$, and then establish the square trees.

\section{The square trees}

In \cite{HW2016-1}, we determined the three cases of squares with kernel $K_m$ (i.e., the maximal kernel word occurring in these squares is $K_m$). 
For $m\geq4$ and $p\geq1$, we denote
$$\Lambda(m,p)=pt_{m-1}+|\mathbb{T}[1,p-1]|_a(t_{m-2}+t_{m-3})+|\mathbb{T}[1,p-1]|_bt_{m-2}.$$
By Property 6.1 in \cite{HW2016-3}, we have the position of the last letter of the $p$-th occurrence of $K_m$ that $P(K_m,p)=\Lambda(m,p)+k_m-1$.
Thus we can define three sets for $m\geq4$ and $p\geq1$, which contain the three cases of squares, respectively.
\begin{equation*}
\begin{cases}
\langle1,K_m,p\rangle&=
\{P(\omega\omega,p):Ker(\omega\omega)=K_m,|\omega|=t_{m-1},\omega\omega\prec\mathbb{T}\}\\
&=\{\Lambda(m,p)+t_{m-1},\cdots,\Lambda(m,p)+k_{m+3}-2\};\\
\langle2,K_m,p\rangle&=
\{P(\omega\omega,p):Ker(\omega\omega)=K_m,|\omega|=t_{m-4}+t_{m-3},\omega\omega\prec\mathbb{T}\}\\
&=\{\Lambda(m,p)+t_{m-3}+t_{m-4},\cdots,\Lambda(m,p)+k_{m+2}-2\};\\
\langle3,K_m,p\rangle&=
\{P(\omega\omega,p):Ker(\omega\omega)=K_m,|\omega|=t_{m-4}-k_{m-3},\omega\omega\prec\mathbb{T}\}\\
&=\{\Lambda(m,p)+k_{m}-1,\cdots,\Lambda(m,p)+2t_{m-4}-1\}.
\end{cases}
\end{equation*}

For $m\geq4$ and $p\geq1$, we consider the sets
\begin{equation*}
\begin{cases}
\Gamma_{1,m,p}=\{\Lambda(m,p)+k_{m+2}-1,\cdots,\Lambda(m,p)+k_{m+3}-2\};\\
\Gamma_{2,m,p}=\{\Lambda(m,p)+k_{m+1}-1,\cdots,\Lambda(m,p)+k_{m+2}-2\};\\
\Gamma_{3,m,p}=\{\Lambda(m,p)+k_{m}-1,\cdots,\Lambda(m,p)+k_{m+1}-2\}.
\end{cases}
\end{equation*}
Obviously, $\langle 1,K_m,p\rangle$ (resp. $\langle 2,K_m,p\rangle$, $\langle 3,K_m,p\rangle$) contains the several maximal (resp. maximal, minimal) elements of $\Gamma_{1,m,p}$ (resp. $\Gamma_{2,m,p}$, $\Gamma_{3,m,p}$).
Moreover $\max\Gamma_{2,m,p}+1=\min\Gamma_{1,m,p}$ and
$\max\Gamma_{3,m,p}+1=\min\Gamma_{2,m,p}$.
Using Lemma 6.4 in \cite{HW2016-3},
comparing minimal and maximal elements in these sets below, we have
\begin{equation*}
\begin{cases}
\Gamma_{1,m,p}=\Gamma_{3,m-1,P(a,p)+1}\cup\Gamma_{2,m-1,P(a,p)+1}\cup\Gamma_{1,m-1,P(a,p)+1},
&m\geq4;\\
\Gamma_{2,m,p}=\Gamma_{3,m-2,P(b,p)+1}\cup\Gamma_{2,m-2,P(b,p)+1}\cup\Gamma_{1,m-2,P(b,p)+1},
&m\geq5;\\
\Gamma_{3,m,p}=\Gamma_{3,m-3,P(c,p)+1}\cup\Gamma_{2,m-3,P(c,p)+1}\cup\Gamma_{1,m-3,P(c,p)+1},
&m\geq6.
\end{cases}
\end{equation*}
Thus we establish the recursive relations for any $\Gamma_{1,m,p}$ ($m\geq4$),
$\Gamma_{2,m,p}$ ($m\geq5$) and $\Gamma_{3,m,p}$ ($m\geq6$).
By the relation between $\Gamma_{i,m,p}$ and $\langle i,K_m,p\rangle$, we get the tree structure of the positions of repeated squares in $\mathbb{T}$, called the square trees.
\begin{equation*}
\begin{cases}
\pi_1\langle1,K_m,p\rangle=\langle3,K_{m-1},\hat{a}\rangle\cup\langle2,K_{m-1},\hat{a}\rangle
\cup\langle1,K_{m-1},\hat{a}\rangle,m\geq4;\\
\pi_2\langle2,K_m,p\rangle=\langle3,K_{m-2},\hat{b}\rangle\cup\langle2,K_{m-2},\hat{b}\rangle
\cup\langle1,K_{m-2},\hat{b}\rangle,m\geq5;\\
\pi_3\langle3,K_m,p\rangle=\langle3,K_{m-3},\hat{c}\rangle\cup\langle2,K_{m-3},\hat{c}\rangle
\cup\langle1,K_{m-3},\hat{c}\rangle,m\geq6.
\end{cases}
\end{equation*}
Here we denote $P(\alpha,p)+1$ by $\hat{\alpha}$ for short, $\alpha\in\{a,b,c\}$.

\begin{figure}
\footnotesize
\setlength{\unitlength}{1mm}
\begin{center}
\begin{picture}(140,43)
\put(6,1){51}
\put(6,4){50}
\put(6,7){49}
\put(6,10){48}
\put(1,13){$\langle1,K_6,1\rangle$}
\put(1,0){\line(0,1){16}}
\put(1,0){\line(1,0){13}}
\put(14,0){\line(0,1){16}}
\put(1,16){\line(1,0){13}}
\put(26,1){51}
\put(26,4){50}
\put(21,7){$\langle1,K_5,2\rangle$}
\put(21,0){\line(0,1){10}}
\put(21,0){\line(1,0){13}}
\put(34,0){\line(0,1){10}}
\put(21,10){\line(1,0){13}}
\put(46,22){44}
\put(46,25){43}
\put(41,28){$\langle2,K_5,2\rangle$}
\put(41,21){\line(0,1){10}}
\put(41,21){\line(1,0){13}}
\put(54,21){\line(0,1){10}}
\put(41,31){\line(1,0){13}}
\put(66,34){40}
\put(66,37){39}
\put(61,40){$\langle3,K_5,2\rangle$}
\put(61,33){\line(0,1){10}}
\put(61,33){\line(1,0){13}}
\put(74,33){\line(0,1){10}}
\put(61,43){\line(1,0){13}}
\put(86,1){51}
\put(81,4){$\langle1,K_4,4\rangle$}
\put(81,0){\line(0,1){7}}
\put(81,0){\line(1,0){13}}
\put(94,0){\line(0,1){7}}
\put(81,7){\line(1,0){13}}
\put(106,13){47}
\put(101,16){$\langle2,K_4,4\rangle$}
\put(101,12){\line(0,1){7}}
\put(101,12){\line(1,0){13}}
\put(114,12){\line(0,1){7}}
\put(101,19){\line(1,0){13}}
\put(126,19){45}
\put(121,22){$\langle3,K_4,4\rangle$}
\put(121,18){\line(0,1){7}}
\put(121,18){\line(1,0){13}}
\put(134,18){\line(0,1){7}}
\put(121,25){\line(1,0){13}}
\put(15,10){\vector(1,-1){5}}
\put(15,10){\vector(2,1){25}}
\put(15,10){\line(1,1){25}}
\put(40,35){\vector(1,0){20}}
\put(16,15){$\pi_1$}
\put(38,7){$\pi_1$}
\put(35,4){\vector(1,0){45}}
\put(35,4){\line(5,1){50}}
\put(85,14){\vector(1,0){15}}
\put(35,4){\line(3,1){54}}
\put(89,22){\vector(1,0){31}}
\put(120,0){\footnotesize{\textbf{(a)}}}
\end{picture}
\begin{picture}(140,47)
\put(6,1){71}
\put(6,4){70}
\put(6,7){69}
\put(6,10){68}
\put(6,13){67}
\put(6,16){66}
\put(6,19){65}
\put(6,22){64}
\put(1,25){$\langle2,K_7,1\rangle$}
\put(1,0){\line(0,1){28}}
\put(1,0){\line(1,0){13}}
\put(14,0){\line(0,1){28}}
\put(1,28){\line(1,0){13}}
\put(26,1){71}
\put(26,4){70}
\put(21,7){$\langle1,K_5,3\rangle$}
\put(21,0){\line(0,1){10}}
\put(21,0){\line(1,0){13}}
\put(34,0){\line(0,1){10}}
\put(21,10){\line(1,0){13}}
\put(46,22){64}
\put(46,25){63}
\put(41,28){$\langle2,K_5,3\rangle$}
\put(41,21){\line(0,1){10}}
\put(41,21){\line(1,0){13}}
\put(54,21){\line(0,1){10}}
\put(41,31){\line(1,0){13}}
\put(66,34){60}
\put(66,37){59}
\put(61,40){$\langle3,K_5,3\rangle$}
\put(61,33){\line(0,1){10}}
\put(61,33){\line(1,0){13}}
\put(74,33){\line(0,1){10}}
\put(61,43){\line(1,0){13}}
\put(86,1){71}
\put(81,4){$\langle1,K_4,6\rangle$}
\put(81,0){\line(0,1){7}}
\put(81,0){\line(1,0){13}}
\put(94,0){\line(0,1){7}}
\put(81,7){\line(1,0){13}}
\put(106,13){67}
\put(101,16){$\langle2,K_4,6\rangle$}
\put(101,12){\line(0,1){7}}
\put(101,12){\line(1,0){13}}
\put(114,12){\line(0,1){7}}
\put(101,19){\line(1,0){13}}
\put(126,19){65}
\put(121,22){$\langle3,K_4,6\rangle$}
\put(121,18){\line(0,1){7}}
\put(121,18){\line(1,0){13}}
\put(134,18){\line(0,1){7}}
\put(121,25){\line(1,0){13}}
\put(15,10){\vector(1,-1){5}}
\put(15,10){\vector(2,1){25}}
\put(15,10){\line(1,1){25}}
\put(40,35){\vector(1,0){20}}
\put(16,15){$\pi_2$}
\put(38,7){$\pi_1$}
\put(35,4){\vector(1,0){45}}
\put(35,4){\line(5,1){50}}
\put(85,14){\vector(1,0){15}}
\put(35,4){\line(3,1){54}}
\put(89,22){\vector(1,0){31}}
\put(120,0){\footnotesize{\textbf{(b)}}}
\end{picture}
\begin{picture}(140,47)
\put(5,7){106}
\put(5,10){105}
\put(5,13){104}
\put(5,16){103}
\put(5,19){102}
\put(5,22){101}
\put(5,25){100}
\put(6,28){99}
\put(6,31){98}
\put(6,34){97}
\put(6,37){96}
\put(1.5,40){$\langle3,K_8,1\rangle$}
\put(1,6){\line(0,1){37}}
\put(1,6){\line(1,0){13}}
\put(14,6){\line(0,1){37}}
\put(1,43){\line(1,0){13}}
\put(25,1){108}
\put(25,4){107}
\put(21.5,7){$\langle1,K_5,5\rangle$}
\put(21,0){\line(0,1){10}}
\put(21,0){\line(1,0){13}}
\put(34,0){\line(0,1){10}}
\put(21,10){\line(1,0){13}}
\put(45,22){101}
\put(45,25){100}
\put(41.5,28){$\langle2,K_5,5\rangle$}
\put(41,21){\line(0,1){10}}
\put(41,21){\line(1,0){13}}
\put(54,21){\line(0,1){10}}
\put(41,31){\line(1,0){13}}
\put(66,34){97}
\put(66,37){96}
\put(61.5,40){$\langle3,K_5,5\rangle$}
\put(61,33){\line(0,1){10}}
\put(61,33){\line(1,0){13}}
\put(74,33){\line(0,1){10}}
\put(61,43){\line(1,0){13}}
\put(85,1){108}
\put(81.5,4){$\langle1,K_4,9\rangle$}
\put(81,0){\line(0,1){7}}
\put(81,0){\line(1,0){13}}
\put(94,0){\line(0,1){7}}
\put(81,7){\line(1,0){13}}
\put(105,13){104}
\put(101.5,16){$\langle2,K_4,9\rangle$}
\put(101,12){\line(0,1){7}}
\put(101,12){\line(1,0){13}}
\put(114,12){\line(0,1){7}}
\put(101,19){\line(1,0){13}}
\put(125,19){102}
\put(121.5,22){$\langle3,K_4,9\rangle$}
\put(121,18){\line(0,1){7}}
\put(121,18){\line(1,0){13}}
\put(134,18){\line(0,1){7}}
\put(121,25){\line(1,0){13}}
\put(15,25){\vector(3,1){45}}
\put(15,25){\vector(1,0){25}}
\put(15,25){\vector(1,-4){5}}
\put(18,28){$\pi_3$}
\put(38,7){$\pi_1$}
\put(35,4){\vector(1,0){45}}
\put(35,4){\line(5,1){50}}
\put(85,14){\vector(1,0){15}}
\put(35,4){\line(3,1){54}}
\put(89,22){\vector(1,0){31}}
\put(120,0){\footnotesize{\textbf{(c)}}}
\end{picture}
\end{center}
\normalsize
\caption{(a)-(c) are square trees from root $\langle1,K_6,1\rangle$, $\langle2,K_7,1\rangle$, $\langle3,K_7,1\rangle$, respectively.}
\label{fig:1}
\end{figure}
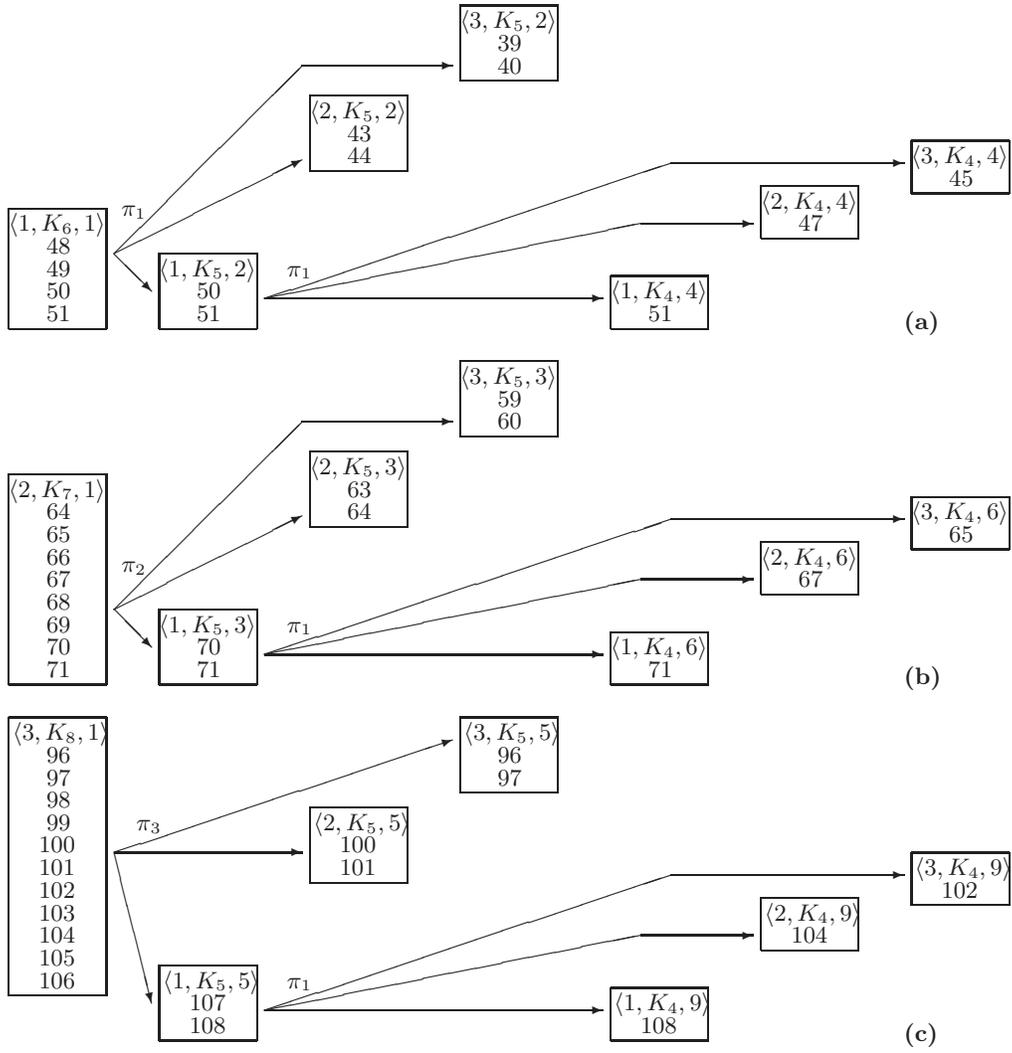

On the other hand,
for $m\geq4$ and $i\in\{1,2,3\}$, each $\langle i,K_{m},1\rangle$ belongs to the square trees. Moreover $\langle i,K_m,\hat{a}\rangle$ (resp. $\langle i,K_m,\hat{b}\rangle$, $\langle i,K_m,\hat{n}\rangle$) is subset of $\pi_1\langle1,K_{m+1},p\rangle$
(resp. $\pi_2\langle2,K_{m+2},p\rangle$, $\pi_3\langle3,K_{m+3},p\rangle$).
Notice that $\mathbb{N}=\{1\}\cup\{P(a,p)+1\}\cup\{P(b,p)+1\}\cup\{P(c,p)+1\}$, the square trees contain all $\langle i,K_m,p\rangle$, i.e. all squares in $\mathbb{T}$.
Fig. \ref{fig:1} shows some examples.

\section{Algorithm: the numbers of repeated squares in $\mathbb{T}[1,n]$}

Denote $a(n)=\sharp\{(\omega,p):\omega_p\omega_{p+1}\triangleright\mathbb{T}[1,n]\}$
the number of squares ending at position $n$.
By Proposition \ref{a}, we can calculate $a(n)$, and obversely calculate $A(n)$ by $A(n)=\sum_{i=4}^n a(i)$.
For $m\geq4$, since $k_m=\frac{t_m-2t_{m-1}+t_{m-2}+1}{2}$,
\begin{equation*}
\begin{cases}
\Gamma_{1,m,1}=\{\frac{t_{m}+2t_{m-1}-t_{m-2}-1}{2},\cdots,
\frac{t_{m}+2t_{m-1}+t_{m-2}-3}{2}\};\\
\Gamma_{2,m,1}=\{\frac{-t_{m}+4t_{m-1}+t_{m-2}-1}{2},\cdots,
\frac{t_{m}+2t_{m-1}-t_{m-2}-3}{2}\};\\
\Gamma_{3,m,1}=\{\frac{t_{m}+t_{m-2}-1}{2},\cdots,
\frac{-t_{m}+4t_{m-1}+t_{m-2}-3}{2}\}.
\end{cases}
\end{equation*}

\begin{proposition}\label{a}\ $a([8])=[1]$, $a([9,10])=[0,1]$, $a([11,\cdots,14])=[0,0,0,1]$,
$a([15,16])=[1,1]$, $a([17,\cdots,20])=[0,0,1,1]$, $a([28,\cdots,31])=[1,1,1,1]$, 
\begin{equation*}
\begin{cases}
a(\Gamma_{1,m,1})
=[a(\Gamma_{3,m-1,1}),a(\Gamma_{2,m-1,1}),a(\Gamma_{1,m-1,1})]
+[\underbrace{0,\cdots,0}_{t_{m-2}-k_{m}+1},\underbrace{1,\cdots,1}_{k_{m}-1}];\\
a(\Gamma_{2,m,1})
=[a(\Gamma_{3,m-2,1}),a(\Gamma_{2,m-2,1}),a(\Gamma_{1,m-2,1})]
+[\underbrace{0,\cdots,0}_{t_{m-3}-k_{m}+1},\underbrace{1,\cdots,1}_{k_{m}-1}];\\
a(\Gamma_{3,m,1})
=[a(\Gamma_{3,m-3,1}),a(\Gamma_{2,m-3,1}),a(\Gamma_{1,m-3,1})]
+[\underbrace{1,\cdots,1}_{t_{m-4}-k_{m-3}+1},\underbrace{0,\cdots,0}_{k_{m-3}-1}].
\end{cases}
\end{equation*}
\end{proposition}

Denote $\Phi_m=\sum a(\Gamma_{3,m,1})+\sum a(\Gamma_{2,m,1})+\sum a(\Gamma_{1,m,1})$.
The immediately corollaries of Proposition \ref{a} are
$\sum a(\Gamma_{1,m,1})=\Phi_{m-1}+k_m-1$,
$\sum a(\Gamma_{2,m,1})=\Phi_{m-2}+k_m-1$,
$\sum a(\Gamma_{3,m,1})=\Phi_{m-3}+t_{m-4}-k_{m-3}+1$.
Moreover, for $m\geq7$,
$$\Phi_m=\Phi_{m-1}+\Phi_{m-2}+\Phi_{m-3}+\tfrac{-3t_{m}+6t_{m-1}+t_{m-2}-1}{2}.$$
By induction, we can prove the 4 properties in Fig \ref{fig:3}.

\begin{figure}
(a) $\Phi_m=\tfrac{m}{22}(-5t_{m}+14t_{m-1}+4t_{m-2})
+\tfrac{1}{44}(67t_{m}-166t_{m-1}+5t_{m-2})+\tfrac{1}{4}$.

(b) $\sum a(\Gamma_{i,4,1})=1$ where $i\in\{1,2,3\}$, and
\begin{equation*}
\begin{cases}
\sum a(\Gamma_{1,m,1})
=\frac{m}{22}(4t_{m}-9t_{m-1}+10t_{m-2})
+\frac{1}{44}(19t_{m}+36t_{m-1}-169t_{m-2})-\frac{1}{4};\\
\sum a(\Gamma_{2,m,1})
=\frac{m}{22}(10t_{m}-6t_{m-1}-19t_{m-2})
+\frac{1}{44}(-189t_{m}+156t_{m-1}+331t_{m-2})-\frac{1}{4};\\
\sum a(\Gamma_{3,m,1})
=\frac{m}{22}(-19t_{m}+29t_{m-1}+13t_{m-2})
+\frac{1}{44}(237t_{m}-358t_{m-1}-157t_{m-2})+\frac{3}{4}.
\end{cases}
\end{equation*}
(c) $\sum_{j=4}^{m-1}\Phi_j
=\tfrac{m}{44}(13t_{m}-10t_{m-1}+5t_{m-2})
+\tfrac{2}{11}(-8t_{m}+8t_{m-1}-7t_{m-2})+\tfrac{m}{4}+2$.
\begin{equation*}
\mathrm{(d)}
\begin{cases}
A(\max\Gamma_{3,m,1})=\frac{m}{44}(-25t_{m}+48t_{m-1}+31t_{m-2})
+\frac{1}{44}(173t_{m}-294t_{m-1}-213t_{m-2})
+\frac{m+11}{4};~~\\
A(\max\Gamma_{2,m,1})=\frac{m}{44}(-5t_{m}+36t_{m-1}-7t_{m-2})
+\frac{1}{22}(-8t_{m}-69t_{m-1}+59t_{m-2})
+\frac{m+10}{4};\\
A(\max\Gamma_{1,m,1})=\frac{m}{44}(3t_{m}+18t_{m-1}+13t_{m-2})
+\frac{1}{44}(3t_{m}-102t_{m-1}-51t_{m-2})
+\frac{m+9}{4}.
\end{cases}
\end{equation*}
\caption{These properties can be proved easily by induction, where (a) and (c) hold for $m\geq4$, (b) and (d) hold for $m\geq5$.}
\label{fig:3}
\end{figure}

Obversely we can calculate $A(n)$ by $A(n)=\sum_{i=4}^n a(i)$. But when $n$ large, this method is complicated. Now we turn to give a fast algorithm.
For any $n\geq52$, let $m$ such that
$n\in\Gamma_{3,m,1}\cup\Gamma_{2,m,1}\cup\Gamma_{1,m,1}
=\{\tfrac{t_{m}+t_{m-2}-1}{2},\cdots,\tfrac{t_{m}+2t_{m-1}+t_{m-2}-3}{2}\}.$
We already determine the expression of $A(\max{\Gamma_{i,m,1}})$ for $i\in\{1,2,3\}$, $m\geq5$. In order to calculate $A(n)$, we only need to calculate $\sum_{i=\min{\Gamma_{i,m,1}}}^n a(i)$.

\begin{figure}
\footnotesize
\textbf{(a)} $n\in\Gamma_{3,m,1}=\{\frac{t_{m}+t_{m-2}-1}{2},\cdots,
\frac{-t_{m}+4t_{m-1}+t_{m-2}-3}{2}\}$ for $m\geq7$. Denote
\begin{equation*}
\begin{cases}
\theta_m^1=\min\Gamma_{3,m,1}=\frac{t_{m}+t_{m-2}-1}{2};\\
\theta_m^2=\min\Gamma_{3,m,1}+|\Gamma_{3,m-3,1}|=\frac{-5t_{m}+10t_{m-1}+3t_{m-2}-1}{2};\\
\theta_m^3=\min\Gamma_{3,m,1}+|\Gamma_{3,m-3,1}|+|\Gamma_{2,m-3,1}|
=\frac{-t_{m}+6t_{m-1}-3t_{m-2}-1}{2};\\
\eta_m^1=\min\Gamma_{3,m,1}+t_{m-4}-k_{m-3}+1=-2t_{m}+5t_{m-1}.\\
\theta_m^4=\max\Gamma_{3,m,1}+1=\min\Gamma_{2,m,1}=\frac{-t_{m}+4t_{m-1}+t_{m-2}-1}{2}.
\end{cases}
\end{equation*}
Obviously, $\theta_m^3<\eta_m^1<\theta_m^4$ for $m\geq7$, and
$\min\Gamma_{3,m,1}-\min\Gamma_{3,m-3,1}=t_{m-1}$. By Property \ref{a}, we have:
for $n\geq52$, let $m$ such that $n\in\Gamma_{3,m,1}$, then $m\geq7$ and
\tiny
\begin{equation*}
\begin{split}
&\sum\nolimits_{i=\min\Gamma_{3,m,1}}^na(i)\\
=&\begin{cases}
\sum_{i=\min\Gamma_{3,m-3,1}}^{n-t_{m-1}}a(i)+n-\min\Gamma_{3,m,1}+1,
&\theta_m^1\leq n<\theta_m^2;\\
\sum_{i=\min\Gamma_{2,m-3,1}}^{n-t_{m-1}}a(i)+\sum a(\Gamma_{3,m-3,1})
+n-\min\Gamma_{3,m,1}+1,&\theta_m^2\leq n<\theta_m^3;\\
\sum_{i=\min\Gamma_{1,m-3,1}}^{n-t_{m-1}}a(i)+\sum a(\Gamma_{3,m-3,1})
+\sum a(\Gamma_{2,m-3,1})+n-\min\Gamma_{3,m,1}+1,&\theta_m^3\leq n<\eta_m^1;\\
\sum_{i=\min\Gamma_{1,m-3,1}}^{n-t_{m-1}}a(i)+\sum a(\Gamma_{3,m-3,1})
+\sum a(\Gamma_{2,m-3,1})+\frac{-5t_{m}+10t_{m-1}-t_{m-2}+1}{2},&otherwise.
\end{cases}
\end{split}
\end{equation*}
\footnotesize
\textbf{(b)} $n\in\Gamma_{2,m,1}=\{\frac{-t_{m}+4t_{m-1}+t_{m-2}-1}{2},\cdots,
\frac{t_{m}+2t_{m-1}-t_{m-2}-3}{2}\}$ for $m\geq6$. Denote
\begin{equation*}
\begin{cases}
\theta_m^5=\min\Gamma_{2,m,1}+|\Gamma_{3,m-2,1}|=\frac{3t_{m}-5t_{m-2}-1}{2};\\
\eta_m^2=\min\Gamma_{2,m,1}+t_{m-3}-k_{m}+1=2t_{m-1}-t_{m-2}.\\
\theta_m^6=\min\Gamma_{3,m,1}+|\Gamma_{3,m-3,1}|+|\Gamma_{2,m-3,1}|
=\frac{3t_{m}-2t_{m-1}-t_{m-2}-1}{2};\\
\theta_m^7=\max\Gamma_{2,m,1}+1=\min\Gamma_{1,m,1}=\frac{t_{m}+2t_{m-1}-t_{m-2}-1}{2}.
\end{cases}
\end{equation*}
Obviously, $\theta_m^5<\eta_m^2\leq\theta_m^6$ for $m\geq6$, and $\min\Gamma_{2,m,1}-\min\Gamma_{3,m-2,1}=t_{m-1}$. By Property \ref{a}, we have:
for $n\geq32$, let $m$ such that $n\in\Gamma_{2,m,1}$, then $m\geq6$ and
\scriptsize
\begin{equation*}
\begin{split}
&\sum\nolimits_{i=\min\Gamma_{2,m,1}}^na(i)\\
=&\begin{cases}
\sum_{i=\min\Gamma_{3,m-2,1}}^{n-t_{m-1}}a(i),
&\theta_m^4\leq n<\theta_m^5;\\
\sum_{i=\min\Gamma_{2,m-2,1}}^{n-t_{m-1}}a(i)+\sum a(\Gamma_{3,m-2,1}),
&\theta_m^5\leq n<\eta_m^2;\\
\sum_{i=\min\Gamma_{2,m-2,1}}^{n-t_{m-1}}a(i)+\sum a(\Gamma_{3,m-2,1})
+n-\eta_m^2+1,&\eta_m^2\leq n<\theta_m^6;\\
\sum_{i=\min\Gamma_{1,m-2,1}}^{n-t_{m-1}}a(i)+\sum a(\Gamma_{3,m-2,1})
+\sum a(\Gamma_{2,m-2,1})+n-\eta_m^2+1,&otherwise.
\end{cases}
\end{split}
\end{equation*}
\footnotesize
\textbf{(c)} $n\in\Gamma_{1,m,1}=\{\frac{t_{m}+2t_{m-1}-t_{m-2}-1}{2},\cdots,
\frac{t_{m}+2t_{m-1}+t_{m-2}-3}{2}\}$ for $m\geq5$. Denote
\begin{equation*}
\begin{cases}
\theta_m^8=\min\Gamma_{1,m,1}+|\Gamma_{3,m-1,1}|=\frac{t_{m}+3t_{m-2}-1}{2};\\
\theta_m^9=\min\Gamma_{1,m,1}+|\Gamma_{3,m-1,1}|+|\Gamma_{2,m-1,1}|
=\frac{-t_{m}+4t_{m-1}+3t_{m-2}-1}{2};\\
\eta_m^3=\min\Gamma_{1,m,1}+t_{m-2}-k_{m}+1=2t_{m-1}.\\
\theta_m^{10}=\max\Gamma_{1,m,1}+1=\min\Gamma_{3,m+1,1}=\frac{t_{m}+2t_{m-1}+t_{m-2}-1}{2}.
\end{cases}
\end{equation*}
Obviously, $\theta_m^9<\eta_m^3<\theta_m^{10}$ for $m\geq5$, and $\min\Gamma_{1,m,1}-\min\Gamma_{3,m-1,1}=t_{m-1}$. By Property \ref{a}, we have:
for $n\geq21$, let $m$ such that $n\in\Gamma_{1,m,1}$, then $m\geq5$ and
\scriptsize
\begin{equation*}
\begin{split}
&\sum\nolimits_{i=\min\Gamma_{1,m,1}}^na(i)\\
=&\begin{cases}
\sum_{i=\min\Gamma_{3,m-1,1}}^{n-t_{m-1}}a(i),
&\theta_m^7\leq n<\theta_m^8;\\
\sum_{i=\min\Gamma_{2,m-1,1}}^{n-t_{m-1}}a(i)+\sum a(\Gamma_{3,m-1,1}),
&\theta_m^8\leq n<\theta_m^9;\\
\sum_{i=\min\Gamma_{1,m-1,1}}^{n-t_{m-1}}a(i)+\sum a(\Gamma_{3,m-1,1})
+\sum a(\Gamma_{2,m-1,1}),&\theta_m^9\leq n<\eta_m^3;\\
\sum_{i=\min\Gamma_{1,m-1,1}}^{n-t_{m-1}}a(i)+\sum a(\Gamma_{3,m-1,1})
+\sum a(\Gamma_{2,m-1,1})+n-\eta_m^3+1,&otherwise.
\end{cases}
\end{split}
\end{equation*}
\normalsize
\caption{(a)-(c) show the three cases of recursive relations between $\sum\nolimits_{i=\min\Gamma_{k,m,1}}^na(i)$ and $\sum\nolimits_{i=\min\Gamma_{t,m-k,1}}^na(i)$, where $k,t\in\{1,2,3\}$, respectively.
These relations are derived directly from the square trees (the tree structure of the positions of repeated squares). Using them, we can calculate $\sum\nolimits_{i=\min\Gamma_{k,m,1}}^na(i)$ fast, and give a fast algorithm for $A(n)$.}
\label{fig:2}
\end{figure}

\vspace{0.2cm}

\noindent\textbf{Algorithm.}
Step 1. For $n\leq51$, calculate $\sum_{i=\min{\Gamma_{i,m,1}}}^n a(i)$ by Property \ref{a}.

Step 2. For $n\geq52$, find the $m$ and $i$ such that $n\in\Gamma{i,m,1}$, then $m\geq7$.
We calculate $\sum_{i=\min{\Gamma_{i,m,1}}}^n a(i)$ by the properties in Fig.\ref{fig:2}.

Step 3. Calculate $A(\min{\Gamma_{i,m,1}}-1)$ by the Property (d) in Fig.\ref{fig:3}.

Step 4. $A(n)=A(\min{\Gamma_{i,m,1}}-1)+\sum_{i=\min{\Gamma_{i,m,1}}}^n a(i)$.

\section{Expression: the numbers of repeated squares in $T_m$}

Since $\theta_m^8\leq t_{m}<\theta_m^9$
and $\theta_{m-1}^6\leq t_{m}-t_{m-1}<\theta_{m-1}^7$ for $m\geq7$, see Fig.\ref{fig:2},
$$\begin{array}{rl} &\sum_{i=\min\Gamma_{1,m,1}}^{t_m}a(i)-\sum_{i=\min\Gamma_{1,m-3,1}}^{t_{m-3}}a(i)\\
=&\sum a(\Gamma_{3,m-1,1})+\sum a(\Gamma_{3,m-3,1})+\sum a(\Gamma_{2,m-3,1})
+2t_{m}-2t_{m-1}-3t_{m-2}+1\\
=&\tfrac{m}{22}(-19t_m+29t_{m-1}+13t_{m-2})
+\tfrac{1}{44}(347t_m-622t_{m-1}-47t_{m-2})+\tfrac{9}{4}.
\end{array}$$
For $m\geq7$, by induction, $\sum_{i=\min\Gamma_{1,m,1}}^{t_m}a(i)$ is equal to
$$\tfrac{m}{44}(23t_m-38t_{m-1}-3t_{m-2})
+\tfrac{1}{44}(-65t_m+164t_{m-1}-105t_{m-2})+\tfrac{3m}{4}-\tfrac{9}{4}.$$
Since $\min\Gamma_{1,m,1}-1=\max\Gamma_{2,m,1}$, $A(t_m)=A(\max\Gamma_{2,m,1})+\sum_{i=\min\Gamma_{1,m,1}}^{t_m}a(i)$.
By the properties in Fig.\ref{fig:2}, we can prove Theorem 21 in in H.Mousavi and J.Shallit\cite{MS2014} in a novel way: for $m\geq3$,
$$A(t_{m})=\tfrac{m}{22}(9t_m-t_{m-1}-5t_{m-2})
+\tfrac{1}{44}(-81t_m+26t_{m-1}+13t_{m-2})+m+\tfrac{1}{4}.$$

\vspace{0.5cm}

\noindent\textbf{\Large{Acknowledgments}}

\vspace{0.4cm}

The research is supported by the Grant NSFC No.11431007, No.11271223 and No.11371210.

\end{CJK*}
\end{document}